\documentclass[12pt]{article}
    \usepackage{amsmath, amssymb, amsthm, amsxtra,    }
\def\be{\begin{equation}}
\def\ee{\end{equation}}

\def\vp{\varphi}
\def\arrowk{^\to{\kern -6pt\topsmash k}}
\def\arrowK{^{^\to}{\kern -9pt\topsmash K}}
\def\arrowt{^\to{\kern -6pt\topsmash t}}
\def\arrowr{^\to{\kern-6pt\topsmash r}}
\def\arrowvp{^\to{\kern -8pt\topsmash\vp}}
\def\tk{\tilde{\kern 1 pt\topsmash k}}
\def\barm{\bar{\kern-.2pt\bar m}}
\def\barN{\bar{\kern-1pt\bar N}}
\def\barA{\, \bar{\kern-3pt \bar A}}

\def\be{\begin{equation}}
\def\ee{\end{equation}}
\numberwithin{equation}{section}
\begin{document}
    \theoremstyle{plain}
    \newtheorem{theorem}{Theorem}
    \newtheorem{lemma}{Lemma}
    \newtheorem{proposition}{Proposition}
    \newtheorem{corollary}[theorem]{Corollary}
    %only theorems and corollaries are on the same counter in this scheme
    \theoremstyle{definition}
    \newtheorem*{definition}{Definition}
    \theoremstyle{remark}
    \newtheorem*{remark}{Remark}
\title{ Arithmetic progressions in multiplicative groups of finite fields }
\author{Mei-Chu Chang\footnote{Research partially
financed by the NSF Grants~DMS~1600154.}\\ \texttt{Department of Mathematics}\\
\texttt{University of California, Riverside}\\\texttt{\small
mcc@math.ucr.edu}}

\date{}
\maketitle
\begin{center}{\bf Abstract}\end{center}

{\it Let $G$ be a multiplicative subgroup of the prime field $\mathbb F_p$ of size $|G|> p^{1-\kappa}$ and $r$ an arbitrarily fixed positive integer. Assuming $\kappa=\kappa(r)>0$ and $p$ large enough, it is shown that any proportional subset $A\subset G$ contains non-trivial arithmetic progressions of length $r$. The main ingredient is the Szemer\'{e}di-Green-Tao theorem. }

\bigskip

\noindent {\bf \large Introduction.}

We denote by $\mathbb F_p$ the prime field with $p$ elements and $\mathbb F_p^*$ its multiplicative group. The main result in this paper is the following.

\medskip

\noindent{\bf Theorem 1.} {\it Given $r\in\mathbb Z^+$, there is some $\kappa=\frac 1{r\,2^{r+1}}>0$ such that the following holds. Let $\delta>0$, $p$ a sufficiently large prime and $G< \mathbb F_p^*$ a subgroup of size
$$|G|>p^{1-\kappa}.$$
Then any subset $A\subset G$ satisfying $|A|>\delta |G|$ contains non-trivial $r$-progressions.}

\medskip

The proof is based on the extension of Szemer\'{e}di's theorem for pseudo-random weights due to Green and Tao, which is also a key ingredient in their proof of arithmetic progressions in the primes. (See \cite{GT}.)
In \S.2 we will recall the precise statement of that result and the various underlying concepts.

The next point is that a multiplicative group behaves like a pseudo-random object (for the relevant notaion of pseudo-randomness). The latter fact is established by rather straightforward applications of Weil's theorem for character sums with polynomial argument. As an introductory result, we illustrate its use by proving

\medskip

\noindent{\bf Proposition 2.}  {\it  Let $r\in\mathbb Z_+$ be fixed, $p$ large enough and $G<\mathbb F_p^*$ a multiplicative group of size
\be\label{0.1}|G|> c_rp^{1-\frac 1{2r}}\ee
Then $G$ contains ${}\;c_r \big(\frac{|G|}p\big)^r p |G|$ many non-trivial $r+1$-progressions.}

\medskip

Taking $r=2$, condition \eqref{0.1} becomes
\be\label{0.2}
|G|>cp^{\frac 34}\ee
ensuring $G$ to contain non-trivial triplets $a, a+b, a+2b$ in arithmetic progression. This last result is simple and well-known, but though one could conjecture a condition of the form $|G|>cp^{\frac 12}$ to suffice, still the best available in this direction. (See \cite{AB} for instance.)

\smallskip

Concerning three-term arithmetic progressions in general sets, we recall Sanders' result \cite{S} which provides the strongest form of Roth's theorem to date and, in the setting of subsets of $\mathbb F_q^n$, $q$ fixed, the solution to the cap set problem due to Ellenberg and Gijswijt \cite{EG}. In negative direction, Behrend's lower bound of $r_3(n)$ has been slightly improved by Elkin \cite{E}. (See also \cite{B}, \cite{CLP}, \cite{M}, \cite{R}, \cite{SS}, \cite{Sz}.)

\smallskip

 It is also natural to expect that when $r$ is large, a condition of the type $|G|>p^{1-\epsilon_r}$ with $\epsilon_r\to 0$ as $r\to \infty$ should be necessary for Proposition 2 to hold. We are not able to show that and could only establish the following.

 \medskip

 \noindent{\bf Proposition 3.} {\it There is a function $\eta_r \to 0$ as $r \to \infty$ and arbitrarily large primes $p$ for which there is a subgroup $G<\mathbb F_p^*$ containing no $r$-progressions and
 \be\label{0.3}
 |G|>p^{\frac 12-\eta_r}.\ee}

 \medskip

 The argument is closely related to a construction in \cite{BGKS}. We note that a more satisfactory result would be \eqref{0.3} with exponent $1-\eta_r$ but $\frac 12-\eta_r$ seems the limit of the method.

 Related to additive shifts of multiplicative subgroups of prime fields, we should also mention the paper of Shkredov and Vyugin \cite{SV}, generalizing results of Konyagin, Heath-Brown and Garcia, Voloch. (See \cite{GV}, \cite{HK}, \cite{K}.)

\medskip

\noindent{\bf\large Notations.}
We recall that the notation $U = O(V )$ is equivalent to the inequality $|U| \le c V$ with some constant $c > 0$, while with the notation $U=o(V)$, in the above inequality, the constant $c$ goes to $0$.
%The notation $f\sim g$ means $f=cg$ for some constant $c$.
We denote by ${\rm ht}F(x)$ the {\it height} of the polynomial $F(x)$, which is the $\max$ of the modulus of the coefficients of $F(x)$. For a set $G$, $\;\mathbb I_G$ is the indicator function of $G$. By $\mathbb E(f\,|\; x\in S)$, we mean the average of $f(x)$ over $x\in S$. The constant $c_r$ is a constant depending on $r$ and may vary even within the same context.

\section{Arithmetic progressions in multiplicative groups.}

 In this section we will prove Proposition 2.

 \smallskip

 First we note that the progression $a, a+b, \ldots, a+rb \in G$  is equivalent to that $a\in G$ and $1+a^{-1}b,\ldots, 1+ra^{-1}b\in G$. Hence we will analyze
\be\label{2.1}
\sum_{x\in\mathbb F_p} \mathbb I_G(1+x)\mathbb I_G(1+2x)\ldots \mathbb I_G(1+rx)\ee
Using the representation
\be\label{1.5}
\mathbb I_G=\frac{|G|}{p-1}\sum_{\chi\equiv 1 \;\rm{ on } \;G}\chi,\ee
we write
\be\label{2.2}
\mathbb I_G=\frac{|G|}{p-1}\bigg(\chi_0+\sum_{\substack{\chi\not=\chi_0\\\chi= 1 \;\rm{ on } \;G}}\chi\bigg).\ee
So we write \eqref{2.1} as
\be\label{1}\bigg(\frac{|G|}{p-1}\bigg)^r(p+\mathcal A),
\ee
where
\be\label{2.4}|\mathcal A|\le \bigg(\frac{p-1}{|G|}\bigg)^r\max \bigg|\sum_{x\in\mathbb F_p} \chi_1(1+x)\ldots \chi_r(1+rx)\bigg|
\ee
with $\max$ taken over all $r$-tuples $\chi_1,\ldots, \chi_r$ of multiplicative characters which are $1$ on $G$ and at least one of them non-trivial.

We now bound the sum in \eqref{2.4}. For the $r$-tuple $\chi_1,\ldots, \chi_r$ obtaining the max, let $I=\{s\in[1, r]: \;\chi_s\not=\chi_0\}$. Assume $\mathcal Y$ generates  $\widehat{\mathbb F_p^*}$ and let $\chi=\mathcal Y^{|G|}$. Then $\chi_s=\chi^{j_s}$, where $j_s< \frac{p-1}{|G|}$. Hence
$$\sum_{x\in\mathbb F_p}\prod_{s\in I} \chi_s(1+sx)=\sum_{x\in\mathbb F_p}\mathcal Y(f(x))$$
with
$$f(x)=\prod_{s\in I} (1+sx)^{j_s|G|}.$$
Since $\mathcal Y$ is of order $p-1$ and $f(x)$ is not a $p-1$-power, Weil's theorem implies
\be\label{2.5}
\bigg|\sum_{x\in\mathbb F_p}\mathcal Y(f(x))\bigg|<|I|\sqrt p.\ee
Assume $|G|>c_rp^{1-\frac 1{2r}}$. It follows that \eqref{1} and hence \eqref{2.1} is bounded below by
\be\label{2.6}
\bigg(\frac{|G|}{p-1}\bigg)^rp-r\sqrt p > c_r \bigg(\frac{|G|}{p-1}\bigg)^rp.\ee
Therefore, $G$ contains at least $c_r\big(\frac{|G|}p\big)^r p |G|$ many non-trivial $r+1$-progressions.

\medskip

\noindent
{\bf Remark 1.1.} We note that if $G\subset \mathbb F_p^*$ is a random set, then the expected size of \eqref{2.1} would also be $\big(\frac{|G|}{p-1}\big)^rp$. So the above observation indicates a random behavior of sufficiently large multiplicative group in terms of $r$-progressions. (This point of view will be exploited further in the next section.)

\section{Progressions in large subsets of multiplicative groups.}

An interesting problem is the following.

\smallskip

\noindent
{\it How large can $G\subset \mathbb F_p^*$ be without containing an r-progression? }

\smallskip

 In this section we will prove Theorem 1. We will use the Green-Tao extension of Szemer\'{e}di's theorem for large subsets of pseudo-random sets. (See Theorem 2.2 in \cite{H}.)

\medskip

\noindent
{\bf Theorem GT.} {\it Let $\nu : \mathbb Z_N \to \mathbb R^+$ be a pseudo-random weight, and let $r\in \mathbb Z^+$. Then for any $\delta >0$, there is $c_r(\delta)>0$ satisfying the following property.

\noindent
For any $f: \mathbb Z_N\to \mathbb R$ such that
\be\label{3.2}0\le f(x)\le \nu(x), \forall x \;\; \text{ and } \;\; \mathbb E(f\,|\; \mathbb Z_N)\ge \delta,\ee
we have
\be\label{3.3}
\mathbb E(f(x)f(x+t)\ldots f(x+rt)\,|\; x, t\in \mathbb Z_N)\ge c_r(\delta)-o(1).\ee}
\noindent (Note that here the notation $\mathbb E$ refers to the normalized sum.)

\medskip

In order to apply this result, one will need to verify that under appropriate assumptions, $\mathbb I_G$ for $G\subset \mathbb F_p^*$, satisfies the required pseudo-randomness conditions.

We call that $\nu$ is a pseudo-random weight if $\nu$ satisfies the following two conditions.

\smallskip

\noindent (1). {\it Condition on linear forms.}

\noindent Let $m_0, t$ and $L\in \mathbb Z$ be constants depending on $r$ only. Let $m\le m_0$ be an integer and $\psi_1,\ldots, \psi_m:\mathbb Z_N^t\to\mathbb Z_N$ be functions of the form
\be\label{3.4}
\psi_i({\bf x})=b_i+\sum_{j=1}^t L_{i,j}x_j,\ee
where ${\bf x}=(x_1,\ldots, x_t)$, $b_i\in\mathbb Z$, $|L_{i,j}|\le L$ and the $m$ vectors $(L_{i,j})_{1\le j \le t} \in\mathbb Z^t$ are pairwisely non-collinear.

\noindent
Then
\be\label{3.5}
\mathbb E\big(\nu(\psi_1({\bf x}))\ldots \nu(\psi_m({\bf x})\big)|\;{\bf x}\in \mathbb Z_N^t)=1+o(1).\ee

\smallskip

\noindent (2). {\it Condition of correlations.}

\noindent Let $q_0\in\mathbb Z$ be a constant. Then there exists $\tau:\mathbb Z_N\to \mathbb R^+$ satisfying
\be\label{3.6}
\text{ for all }\;\ell\ge 1, \mathbb E(\tau^{\ell}(x)\,|\; x\in \mathbb Z_N)=O_{\ell}(1)\ee
such that for all $q\le q_0$ and $h_1, \ldots , h_q\in\mathbb Z_N$ (not necessarily distinct), we have
\be\label{3.7}
\mathbb E\big(\nu(x+h_1)\nu(x+h_2)\ldots\nu(x+h_q)\,|\;x\in \mathbb Z_N\big)\le\sum_{1\le i\le j\le q}\tau(h_i-h_j).\ee

\smallskip

\noindent{\bf Remark 2.1.} As Y. Zhao pointed out that in his paper \cite{CFZ} with D. Conlon and J. Fox, they showed that in applying Theorem GT one only needs to verify the $m_0$-linear forms condition (with $m_0=r\,2^{r-1}$), and that the correlation condition is actually unnecessary.

\medskip

\noindent {\bf Proof of Theorem 1.} In our application of Theorem GT, $\mathbb Z_N$ will be $\mathbb F_p$ with additive structure and $\nu=\frac{p-1}{|G|}\mathbb I_G$. We will verify the condition on linear forms above by using Weil's theorem.

 \smallskip

\noindent Using the representation \eqref{2.2}, we have
\be\label{3.8}\begin{aligned}
\nu=\;&\;\frac {p-1}{|G|} \mathbb I_G\;=\; \frac {p-1}{|G|}\;\frac {|G|} {p-1}\sum_{\chi= 1 \;\rm{ on } \;G}\chi\\=\;&
\;\chi_0+\sum_{\substack{\chi\not=\;\chi_0\\\chi= 1 \;\rm{ on } \;G}}\chi.\end{aligned}\ee

In \eqref{3.5}, the trivial character $\chi_0$ contributes for $1$ and the additional contribution may be bounded as in \S 1 by
\be\label{3.9}
\bigg(\frac {p-1}{|G|}\bigg)^m p^{-t} \max\bigg|\sum_{{\bf x}\in\mathbb F_p^t}\chi_1(\psi_1({\bf x}))\ldots \chi_m(\psi_m({\bf x}))\bigg|
\ee
with $\max$ taken over all $m$-tuples $\chi_1,\ldots, \chi_m$, which are $1$ on $G$ and not all $\chi_0$. For the $m$-tuples $\chi_1,\ldots, \chi_m$ obtaining the max, let $I=\{s\in[1,m]:\chi_s\not=\chi_0\}$, hence $\chi_s=\mathcal Y^{j_s|G|}$, with $j_s< \frac{p-1}{|G|}$ for $s\in I$. We obtain
\be\label{3.10}
\sum_{{\bf x}\in\mathbb F_p^t}\chi_1(\psi_1({\bf x}))\ldots \chi_m(\psi_m({\bf x}))=\sum_{{\bf x}\in\mathbb F_p^t}\mathcal Y\big(\prod_{s\in I}\psi_s({\bf x})^{j_s|G|}\big)\ee
To introduce a new variable $z$, we perform a shift ${\bf x}\mapsto {\bf x} + z{\bf a}$, where ${\bf a}\in \{1,\ldots, m\}^t$ may be chosen such that
\be\label{3.11}
\sum_{j=1}^t L_{s,j}a_j\not= 0,\; \text{ for } s=1, \ldots, m.\ee
Recall that $|L_{s,j}|\le L$ and the m vectors $(L_{s,j})_{j=1,\ldots, t}\in\mathbb Z^t$ are pairwisely non-collinear. Hence we may choose ${\bf a}$ as above to fulfill \eqref{3.11} and moreover $\sum_{j=1}^t L_{s,j}a_j\not\equiv 0\pmod p$. We estimate \eqref{3.10} as
\be\label{3.12}
\frac 1p\sum_{{\bf x}\in\mathbb F_p^t}\bigg|\sum_{z=0}^{p-1}\mathcal Y(f_{\bf x} (z))\bigg|,\ee
where
\be\label{3.13}
f_{{\bf x}} (z)=\prod_{s\in I}\big(\big(\sum_jL_{s,j}a_j\big)z+\psi_s({\bf x})\big)^{j_s|G|}.\ee
Clearly, $f_{{\bf x}} (z)$ will not be a ($p-1$)-power of a polynomial, if the following expressions
\be\label{3.14}
\frac {\psi_s({\bf x})}{\sum_jL_{s,j}a_j}, \;\; s\in I\ee
are pairwisely distinct.

To estimate the double sum in \eqref{3.12}, we write $\sum_{{\bf x}\in\mathbb F_p^t}$ as $\sum^{(1)}+\sum^{(2)}$, where $\sum^{(1)}$ is over those ${\bf x}\in \mathbb F_p^t$ for which \eqref{3.14} are pairwisely distinct and $\sum^{(2)}$ over the other ${\bf x}$.

By Weil's theorem
\be\label{3.15}
\frac 1p\sum{}^{^{(1)}}\bigg|\sum_{z=0}^{p-1}\mathcal Y(f_{\bf x} (z))\bigg|\le |I|\; p^{t-1}\sqrt p.\ee

For $\sum^{(2)}$ we estimate trivially.
\be\label{3.16}\begin{aligned}&
\frac 1p\sum{}^{^{(2)}}\bigg|\sum_{z=0}^{p-1}\mathcal Y(f_{\bf x} (z))\bigg|\\\le\;\; &\sum_{\substack{s, s'\in I\\s\not= s'}}\bigg|\bigg\{{\bf x}\in \mathbb F_p^t:\frac {\psi_s({\bf x})}{\sum_jL_{s,j}a_j}=\frac {\psi_{s'}({\bf x})}{\sum_jL_{s',j}a_j}\bigg\}\bigg|\end{aligned}\ee
Since $(L_{s,j})_{1\le j\le t}$ and $(L_{s',j})_{1\le j\le t}$ are not collinear (and bounded), there is some $j_0$ such that
$$\frac{L_{s, j_0}}{\sum L_{s,j}a_j}\;-\;\frac{L_{s', j_0}}{\sum L_{s',j}a_j}\in\mathbb F_p^*\;.$$
This shows that \eqref{3.16} is bounded by $r^2p^{t-1}$. Therefore, we proved that \eqref{3.12} is bounded by $rp^{t-\frac 12}$, and \eqref{3.9} is bounded by
\be\label{3.18}
c_r \frac {(p-1)^m}{|G|^m \sqrt p},\ee
which is bounded by $p^{-\frac 14}$, assuming
\be\label{3.19} |G|> p^{1-\frac 1{4m_0}}.\quad\square\ee
\section{Construction of large multiplicative groups with no $r$-progressions.}

In this section we will prove Proposition 3. Our argument is very similar to the proof of Theorem 39 in \cite{BGKS}, where it is shown that there is a subset $\Delta\subset \mathcal P_T=\{p: p \text{ is a prime, and } p\le T\}$, $|\Delta |<\delta\frac T{\log T}$ with $\delta=\delta(r)\to 0$ as $r\to\infty$ and such that for any $p\in\mathcal P_T\setminus \Delta$ and any $t\in\mathbb Z$
\be\label{4.3}
\max \big({\rm ord}_p(t+1),\ldots, {\rm ord}_p(t+r)\big)>T^{\frac 12-\delta}.\ee
Obviously, \eqref{4.3} implies that
\be\label{4.4}
{\rm ord}_p \langle t+1,\ldots, t+r\rangle > T^{\frac 12-\delta},\ee
which is the only relevant property for us.

As in \S1, if $a, a+b, \ldots, a+rb\in G\subset\mathbb F_p^*,$ and $\, b\in F_p^*$, then $1+t, 1+2t, \ldots, 1+rt\in G, t\equiv a^{-1}b\pmod p$ and hence we obtain $t\in\mathbb Z, t\not\equiv 0\pmod p$
 such that
$${\rm ord}_p\langle 1+t,\ldots,1+rt\rangle \le |G|.$$
Thus our purpose is to ensure that for all $t\not\equiv 0\pmod p$ such that
\be\label{4.5}
{\rm ord}_p\langle 1+t,\ldots,1+rt\rangle >p^{\frac 12-\delta},\ee with $p$ such that $p-1$ has a divisor $d$ in the interval $[p^{\frac 12-\eta}, p^{\frac 12-\delta}]$.  Then the subgroup $G<\mathbb F_p^*$ of order $d$ will have no ($r+1$)-progression. Assuming \eqref{4.5} holds for all $p\in\mathcal P_T\setminus \Delta$ with $|\Delta|<\delta\frac T{\log T}$, it will then suffice (taking $\eta=c\delta$) to invoke

\medskip

\noindent{\bf Lemma 3.1.} {\it Let notations be as above. Then
\be\label{4.6}
\big|\big\{ p\in\mathcal P_T: p-1 \text{ has a prime divisor in the interval } [T^{\frac 12-\eta}, T^{\frac 12-\frac {\eta}2}]\big\}\big|> c\eta\frac T{\log T}.\ee}

\smallskip

\noindent{\bf Proof.} In Bombieri-Vinogradov theorem, taking
\be\label{4.7}
Q=T^{\frac 12}\big(\log T\big)^{-10}\ee
we have
\be\label{4.8}
\sum_{q\le Q}\bigg|\psi(T; q, 1)-\frac T{\phi(q)}\bigg|=O\big(T^{\frac 12} Q \big(\log T\big)^5\big)<cT\big(\log T\big)^{-5},\ee where $\phi(q)$ is the Euler's totient function and$$\psi(T; q, 1)=\sum_{\substack{n\leq T\\ n\equiv 1 \mod q}}\Lambda(n),$$
$\Lambda (n)$ being the von Mangoldt function.
Denote
$$\Omega=\bigg\{q\in \big[T^{\frac 12-\eta}, T^{\frac 12-\frac {\eta}2}\big]\cap \mathcal P:\psi(T; q, 1)<\frac T{2\phi(q)}\bigg\}.$$
Let $[2^k, 2^{k+1}]\subset [T^{\frac 12-\eta}, T^{\frac 12-\frac {\eta}2}]:=I$. From \eqref{4.8},
$$\big|\Omega \cap [2^k, 2^{k+1}]
\big|\;\frac T{2^{k+1}}< c\;T\big(\log T\big)^{-5},$$
hence
\be\label{4.9}
\big|\Omega \cap [2^k, 2^{k+1}]
\big|\;< c\frac {2^k}{\big(\log T\big)^5} <\frac 1{100}\big|\mathcal P\cap [2^k, 2^{k+1}]\big|.\ee

Clearly, \eqref{4.9} and the prime number theorem imply that
$$\begin{aligned}
\sum_{\substack{q\not\in\Omega\\q\in I\cap \mathcal P}}\,\frac 1q\; =\;\;&\;\sum_{(\frac 12-\eta)\log T<k<(\frac 12-\frac {\eta}2)\log T}\;\;\sum_{\substack{q\not\in\Omega\\q\in [2^k, 2^{k+1}]\cap \mathcal P}}\,\frac 1q\quad\qquad\\
<\;\;&\;\sum_{(\frac 12-\eta)\log T<k<(\frac 12-\frac {\eta}2)\log T}\;\frac 1{2^k}\;\; \big|\mathcal P\cap [2^k, 2^{k+1}]\big|\;<2\eta .\qquad\qquad\end{aligned}$$
Let $\sigma <2 \eta$ be a parameter (to be specified). From the preceding, there is a subset $S\subset I\cap \mathcal P$, $S\cap \Omega=\emptyset$, such that
\be\label{4.10}
\sigma<\sum_{q\in S}\frac 1q< 2\sigma,\ee
and since $S\cap\Omega=\emptyset$, we have for all $q\in S$
\be\label{4.11}
|A_q|\ge \frac T{2\big(\log T\big)q}, \;\;\text{where } A_q:=\{p<T: p\equiv 1\pmod q\}.\ee
From the inclusion/exclusion principle and the Brun-Titchmarsh theorem, the left hand side of \eqref{4.6} is at least
\be\label{3}\begin{aligned}
\bigg|\bigcup_{q\in S}A_q\bigg|\ge &\; \sum_{q\in S}|A_q|-\sum_{\substack{q_1, q_2\in S\\q_1\not= q_2}}|A_{q_1q_2}|\\
\ge &\; \frac T{2\log T}\sum_{q\in S}\frac 1q-\sum_{\substack{q_1, q_2\in S\\q_1\not= q_2}}\bigg\{\frac {2T}{\phi(q_1 q_2)\log \frac T{q_1q_2}}\;\bigg(1+O\bigg(\frac 1{\log\frac T{q_1q_2}}\bigg)\bigg)\bigg\}\end{aligned}\ee
Since $\phi(q_1q_2)=(q_1-1)(q_2-1),$ and $ q_1q_2\le T^{1-\eta}$ for $q_1\not= q_2$ in $S$, \eqref{3} is bounded below by
$$\begin{aligned} &\;\;\frac T{\log T}\;\bigg(\frac 12\sum_{q\in S} \frac 1q-\frac 3{\eta}\bigg(\sum_{q\in S}\frac 1q\bigg)^2\bigg)\qquad\qquad\qquad\qquad\qquad\\
=&\;\;\frac T{\log T}\;\bigg(\frac{\sigma}2-\frac 3{\eta}\sigma^2\bigg)\,>c\,\eta\frac T{\log T}\end{aligned}
$$
for some small $c>0$ and appropriate choice of $\sigma.\quad\square$

\medskip

Returning to the proof of Theorem 39 in \cite{BGKS}, a key ingredient is Lemma 17 (in \cite{BGKS}) depending on a result from \cite{ESS} on additive relations in multiplicative subgroups of $\mathbb C^*$. Keeping \eqref{4.5} in mind, the appropriate variant of Lemma 17 we will need is the following.

\medskip

\noindent{\bf Lemma 3.2.} {\it Let $z\in \mathbb C^*$ and $r\in \mathbb Z_+$ be sufficiently large. Consider the set $\mathcal A=\{1+sz: 1\le s\le r\}\subset \mathbb C$. Then there is a multiplicative independent subset $\mathcal A_0\subset \mathcal A$ of size
\be\label{4.13}
|\mathcal A_0|>\,c\,\log r.\ee}

\medskip

The proof is the same as Lemma 17 in \cite{BGKS}. Note that one distinction is that we have to assume $z \not= 0$, which will also lead to a small modification in the proof of Theorem 39 in \cite{BGKS}, in order to establish \eqref{4.5}. Thus

\medskip

\noindent{\bf Lemma 3.3.} {\it There is a subset $\Delta\subset \mathcal P_T, \, |\Delta|= o\big(\frac T{\log T}\big)$ such that every $p\in \mathcal P_T \setminus \Delta$ has the following property.

\noindent If $t\in\mathbb Z, t\not\equiv 0\pmod p$, then
\be\label{4.15}
{\rm ord}_p\langle 1+t, \ldots, 1+rt\rangle >p^{\frac 12-\delta},\ee
where $\delta=\delta(r)$.}
% as $r\to \infty$.}

\medskip

\noindent{\bf Proof.} The basic strategy is the same as that of Theorem 39 in \cite{BGKS}.

 We fix an integer $r_0=[\log r]$,  let \be\label{4.16}
\delta=\delta(r)=\frac{100}{r_0},\ee
and choose $u\in\mathbb Z_+$ such that
\be\label{4.17}
u^{r_0}=c T^{\frac 12-\delta} \;\text{ and }\;\; \frac 12T^{\frac 12-\delta}<u^{r_0}<2  T^{\frac 12-\delta} .\ee

Let $\mathcal E$ be the collection of all subsets $E\subset\{1,\ldots,r\}, |E|=r_0$.

Next, given any two subsets $E_1, E_2\subset\{1,\ldots r\},\, E_1\cap E_2=\emptyset, \; 0<|E_1|+|E_2|\le r_0$, and exponents $\tilde{u}=(u_s)_{s\in E_1\cup E_2}, 1\le u_s\le u$, we introduce the polynomial
\be\label{4.18}
F=F_{E_1, E_2,\tilde{u}}(x)=\prod_{s\in E_1}(1+sx)^{u_s}-\prod_{s\in E_2}(1+sx)^{u_s}\in\mathbb Z[x].\ee

Note that $x$ is always a factor of $F(x)$. Clearly $\deg F(x)\le r_0u$, ${\rm ht} F(x)\le r^{\,2r_0u}$, and there are at most $2^{r_0}\binom{r\;}{r_0}u^{r_0}$ such polynomials.

Denote by $\mathcal F\subset\mathbb Z[x]$ the collection of all irreducible factors $f(x)\in\mathbb Z[x]$ and $f(x)\not=x$ extracted from all polynomials of the form \eqref{4.18}. Hence
\be\label{4.19}
|\mathcal F|\le r_02^{r_0}\binom{r}{r_0}u^{r_0+1}.\ee

Next, if $f,g\in\mathcal F, f\not\sim g$, (i.e. $f$ and $g$ are not proportional) then the resultant of $f, g$ satisfies
\be\label{4.20}
{\rm Res}(f,g)\in \mathbb Z\setminus \{0\} \;\;\text{ and }\;\; |{\rm Res}(f, g)|<r^{2(r_0u)^2}.\ee
From \eqref{4.19}
\be\label{4.21}
B=\prod_{\substack{f,g\in \mathcal F\\f\not\sim g}}{\rm Res}(f,g)\;\in\mathbb Z\setminus\{0\}\ee
satisfies
\be\label{4.22}
|B|<r^{\;2r_0^4\;4^{r_0}\binom{r\;}{r_0}^2u^{2r_0+4}}<r^{\;r^{2r_
0}\;u^{2r_0+4}}.\ee
By \eqref{4.17} and \eqref{4.16}, for $T$ sufficiently large,  we can bound the exponent in \eqref{4.22} as
$$r^{2r_0}\;u^{2r_0+4}<r^{2r_0}T^{1-2\delta+\frac 2{r_0}}<T^{1-\delta}=o\bigg(\frac T{\log T}\bigg).$$

Therefore, there is a set $\Delta\subset\mathcal P_T$ of primes $p\le T,$ with $ |\Delta|=o\big(\frac T{\log T}\big)$ such that $(p, B)=1$ for all $p\in \mathcal P_T\setminus\Delta$.

Now, take $p\in \mathcal P_T\setminus\Delta$ and suppose there exists some $t\in\mathbb Z, \; t\not\equiv 0\pmod p$ such that
$${\rm ord}_p\langle 1+t, \ldots, 1+rt\rangle <u^{r_0}.$$
Then, for all $E\in\mathcal E$, there are $E_1, E_2\subset E,\; E_1\cap E_2=\emptyset, \; |E_1|+|E_2|\ge 1$ and $\tilde{u}=(u_s)_{s\in E_1\cup E_2}$ such that $F_{E_1, E_2,\tilde{u}}(t)\equiv 0\pmod p$. Hence there is a factor $f_E(x)$ of $F_{E_1, E_2,\tilde{u}}(x)$ such that $f_E(t)\equiv 0\pmod p$. Since $t\not\equiv 0\pmod p$, $f_E(x)\not=x$. For all $E, F\in\mathcal E$, since $f_E(x), f_F(x)$ have common root $t\pmod p$
\be\label{4.23}
{\rm Res}(f_E, f_F)\equiv 0\pmod p.\ee
If $f_E\not= c f_F$, then ${\rm Res}(f_E, f_F)|B$, contradicting $(B, p)=1$. Thus $f_E=c f_F$ for all $E, F\in \mathcal E$ and hence have a common root $z\in\mathcal C^*$. But by Lemma 3.2, there is a set $E\in\mathcal E$ such that $\{1+sz:\, s\in E\}$ are multiplicatively independent, implying $F_{E_1, E_2,\tilde{u}}(z)\not=0, f_E(z)\not=0$, which is a contradiction. $\qquad\square$
\section{}
\bigskip
\noindent{\bf Ackonwledgement.} The author would like to thank Yufei Zhao for bringing her attention to \cite{CFZ}. The author would also like to thank the referees for careful reading, which improved an earlier version of the paper.


\begin{thebibliography}{99}

\bibitem{AB}N. Alon, J. Bourgain,  \textit{Additive Patterns in Multiplicative Subgroups}, Geom. Funct. Anal. 24(3), 721-739, (2014).

\medskip

\bibitem{BK}M. Bateman and N. Katz, \textit{New bounds on cap sets}, J. Amer. Math. Soc. 25(2), 585-613,(2012).

\medskip

\bibitem{B}F. A. Behrend, \textit{On sets of integers which contain no three terms in arithmetical progression}, Proc. Natl. Acad. Sci. USA, 32(12), 331-332, (1946).

    \medskip

\bibitem{BGKS}J. Bourgain,M. Z. Garaev, S. V. Konyagin, I. Shparlinski,  \textit{Multiplicative congruences with variables from short intervals}, J. Anal. Math. 124(1), 117-147, (2014).

\medskip

\bibitem{CFZ}D. Conlon, J. Fox, and Y. Zhao, \textit{A relative Szemerédi theorem}, Geom. Funct. Anal. 25, 733–762, (2015).

\medskip

\bibitem{CLP} E. Croot, V. Lev, and P. P. Pach, \textit{Progression-free sets in $\mathbb Z_n^4$ are exponentially small}, preprint, (2016). arXiv:1605.01506.

\medskip

\bibitem{ESS} J.-H. Evertse, H. Schlickewei, W. Schmidt, \textit{Linear equations in variables which lie in a
multiplicative group}, Ann. of Math. (2), 155, 807-836, (2002).


\medskip

\bibitem{GT} B. Green, T. Tao, \textit{The primes contain arbitrarily long arithmetic progressions}, Ann. of Math. (2), 167(2), 481-547, (2008).
\medskip

\bibitem{H} B. Host, \textit{Arithmetic progressions in primes}, S\'{e}minaire Bourbaki 47,  229-246 (2004/2005).

\medskip

%\bibitem{IK}H. Iwaniec, E. Kowalski, \textit{Analytic number theory}, Amer.  Math.  Soc., Providence, RI, (2004).

\medskip

\bibitem{E}M. Elkin,  \textit{ An improved construction of progression free sets}, Israel J. Math. 184, 93-128, (2011).

\medskip

\bibitem{EG} J. Ellenberg, D. Gijswijt, \textit{On large subsets of $\mathbb F_q^n$ with no three-term arithmetic progression}, preprint, (2016). arXiv:1605.09223.

\medskip
%P. Erd˝os and P. Tur´an. On some sequences of integers. J. London Math. Society, 11:261–264,1936

\bibitem{FGR}P. Frankl, R. L. Graham, V. R$\ddot{o}$dl, \textit{On Subsets of Abelian Groups with No 3-Term
Arithmetic Progression}, J. Combin. Theory Ser. A 45, 157-161, (1987).

\medskip

\bibitem{GV}A. Garcia, J.F. Voloch, \textit{Fermat curves over finite fields}, J. Number Theory 30, 345-356, (1988).

\medskip

\bibitem{HK}D. R. Heath-Brown, S. Konyagin, \textit{New bounds for Gauss sums derived from kth powers,
and for Heilbronn’s exponential sum}, Quart. J. Math. 51, 221-235, (2000).

\medskip

\bibitem{K}S. V. Konyagin, \textit{Estimates for trigonometric sums and for Gaussian sums}, IV International
conference ”Modern problems of number theory and its applications”. Part 3, 86-114, (2002).

\medskip

\bibitem{M}R. Meshulam, \textit{On subsets of finite abelian groups with no 3-term arithmetic progressions}, J. Combin. Theory Ser. A 71(1), 168-172, (1995).

\medskip

\bibitem{R}K. Roth, \textit{On certain sets of integers}, J. Lond. Math. Soc. (2), 28, 245-252, (1953).

\medskip

\bibitem{SS}R. Salem, D. Spencer, \textit{On sets of integers which contain no three in arithmetic progression}, Proc. Natl. Acad. Sci. USA, 28, 561-563, (1942).

\medskip

\bibitem{S}T. Sanders, \textit{On Roth's theorem on progressions},  Ann. of Math. (2), 174 (1), 619-636, (2011).

\medskip

\bibitem{SV}I. Shkredov, I. Vyugin, \textit{On additive shifts of multiplicative subgroups}, Mat. Sb. 203(6), 81-100, (2012).

\medskip

\bibitem{Sz}E. Szemer$\acute{e}$di, \textit{On sets of integers containing no k elements in arithmetic progression}, Acta Arith. 27, 299-345, (1975).


    \end{thebibliography}
\end{document}